\documentclass[preprint,12pt]{elsarticle}

\usepackage{amssymb}
\usepackage{algpseudocode}
\usepackage{algorithm}
\usepackage{bm}
\usepackage{listings}

\usepackage{lineno}
\usepackage{pgfplots}
\usepackage{pstricks,pst-node, pst-fun}
\usepackage{dsfont}
\usepackage{amssymb}
\usepackage[british,UKenglish]{babel}
\usepackage[utf8]{inputenc}
\usepackage[a4paper,left=3cm,right=2cm,top=2.5cm,bottom=2.5cm]{geometry}
\usepackage{amsmath}
\usepackage{relsize}
\usepackage{setspace}
\usepackage{subfig}

\DeclareFixedFont{\ttb}{T1}{txtt}{bx}{n}{12} 
\DeclareFixedFont{\ttm}{T1}{txtt}{m}{n}{12}  

\usepackage{color}
\definecolor{deepblue}{rgb}{0,0,0.5}
\definecolor{deepred}{rgb}{0.6,0,0}
\definecolor{deepgreen}{rgb}{0,0.5,0}

\usepackage{listings}

\newcommand\pythonstyle{\lstset{
language=Python,
basicstyle=\ttm,
morekeywords={self,False,True},              
keywordstyle=\ttb\color{deepblue},
emph={MyClass,__init__},          
emphstyle=\ttb\color{deepred},    
stringstyle=\color{deepgreen},
frame=tb,                         
showstringspaces=false
}}

\lstnewenvironment{python}[1][]
{
\pythonstyle
\lstset{#1}
}
{}

\newcommand\pythoninline[1]{{\pythonstyle\lstinline!#1!}}

\journal{Journal Name}

\begin{document}

\begin{frontmatter}





\title{Mechanostat-type effective density correction for Carter-Hayes growth: application to topology optimization and its efficient interpolation for a target strain energy and volume fraction}


\author{Luis Irastorza-Valera$^{a}$, Ricardo Larraínzar-Garijo$^{b,c}$, Javier Montoya-Adárraga$^{b}$, Luis Saucedo-Mora$^{a}$}

\address{$^a$ E.T.S. de Ingeniería Aeronáutica y del Espacio, Universidad Politécnica de Madrid, Pza. Cardenal Cisneros 3, 28040 Madrid, Spain}
\address{$^b$ Orthopaedics and Trauma Department, Hospital Universitario Infanta Leonor, C/ Gran Vía del Este 80, 28031 Madrid, Spain}
\address{$^c$ Surgery Department, Faculty of Medicine, Universidad Complutense de Madrid, Pl. de Ramón y Cajal s/n, 28040 Madrid, Spain}

\begin{abstract}
In automotive and aerospace industries, the need for optimized structures offering the best mechanical performance for the minimum weight is ubiquitous. To that aim, Topology optimization (TO) is a very popular structural design tool. Particularly, the Solid Isotropic Material with Penalization (SIMP) method offers a trade-off between minimum compliance (i.e., maximum stiffness) and a fixed amount material for a given set of static, deterministic boundary conditions. Since TO is a non-convex problem, its gradient can be tuned by filtering the topology's contour, creating sharper material profiles without necessarily compromising optimality. However, despite simplifying the layout, some filters fail to address manufacturability concerns such as capillarity (thin tweaks as struts) generated by uncertain loading, vibration or fatigue.

A tailored density-based filtering strategy is offered to tackle this issue. Additionally, volume fraction is left unconstrained so material can be strategically replenished through a logarithmic rule acting on the updated compliance. In doing so, an interpolation space with three degrees of freedom (volume, compliance, minimum thickness) is created, yielding diverse topologies for the same boundary conditions and design values along different stages of evolving topological families with distinct features. 

The optimization process is further accelerated by introducing the volume-compliance iterative scheme as a physical loss function in a Double Distance Neural Network (D$^2$NN), obtaining similar results to 2,000 steps worth of vanilla iteration within 500 training epochs.
This proposal offers a novel topology optimization design space based on minimum strut thickness - via filtering - and topological families defined by minimum volume fraction and compliance. The methodology is tested on several examples with diverse loading and boundary conditions, obtaining similarly satisfactory results, and then boosted via Machine Learning, acting as a fast and cheap surrogate.

\end{abstract}


\end{frontmatter}

\section{Introduction}

Before the Industrial Revolution, vernacular structural design (housing and bridges) offered utilitarian feasible solutions where simplicity and resistance were equally valued. These designs often failed due to then unknown or under-studied phenomena such as buckling and fatigue, and were mostly oversized for their intended task. 

Nowadays, structures for industrial applications must be sturdy enough to withstand service specifications (technical requirement) while using the minimum necessary material (economic constraint), reaching a compromise between both criteria. The search for the optimum is now mathematically defined, depending on measurable physical variables.

After briefly introducing the state of the art in this first section, the novelties will be presented alongside several illustrative examples: material density correction in Section 2, the proposed erosion filter in Section 3, the growth-based optimization approach in Section 4 and a particular application for vibration purposes in Section 5. Section 6 will cover some relevant case studies applying the suggested methodology. Section 7 offers a Machine Learning surrogate to accelerate the optimization process and Section 8 will offer some conclusions and proposals for future research lines.

\subsection{Topology Optimization}

Topology optimization (TO) is a widely established technique \cite{Michell1904,Rozvany1972,Bendse1988,Bendse1989} for structural design problems where mechanical requirements are subjected to material restrictions by creating a topological object (adding holes) from an initial isotropic bulk block, yielding a binary distribution (material/void). 

Among different techniques, the Solid Isotropic Material with Penalization (SIMP) \cite{Rozvany1972,Bendse1988} is arguably the most common, due to its simplicity. It consists of an iterative minimization of compliance $\mathbf{c=u^T Ku}$ (twice the strain energy $\Psi$, i.e. stiffness maximization) subjected to a target, fixed volume fraction $v_f$ via Young's modulus decrease $\mathbf{E_{t+1}= E_0 \rho_t^p}$, where $\rho$ is the material's density (from 0 - void - to 1 - full -) and $p$ is a penalization parameter. 

Alas, most topology optimization methods (SIMP, ESO, OMP, etc.) \cite{Rozvany2008} can produce theoretically optimized structures which are not easily manufactured or feasible, especially for additive manufacturing (AM): discontinuities, acute angles, fragile struts, etc. This compromises durability as well, since thin struts are more prone to deteriorate and ultimately break. Plus, they could buckle if they are slender enough. Often times, they appear as a result of load or boundary condition variability. On top of that, they do not actually contribute significantly to supporting the load, conversely adding to the compliance sought to be minimized.

To tackle these shortcomings, many filters address material density \cite{Bourdin2001,Sigmund2007filter,Lazarov2010,Andreassen2010}. They can avoid numerical limitations like checkerboard patterns and mesh-dependence and practical ones like resolution and continuity, with various penalization schemes \cite{Liu2016}. Some of them revolve around erosion phenomena, i.e., scratching material off the structure's contour to get a clear void-material frontier, as in nature \cite{Ostanin2017}. Their effects are usually global, targeting all the contour simultaneously instead of specific areas of interest along the iterative process (a sort of evolutionary design). This kind of density filters do in fact enhance manufacturability as they sharpen the void-material boundaries, although they still do not selectively target undesired ribs - a consequence of load variability \cite{IrastorzaValera2025prob}. 

Many tools have been suggested with manufacturability in mind, particularly after the rise of AM and its implicit multi-scalar challenges \cite{Lazarov2016} (length, connectivity, self-support) and opportunities (e.g., tailored anisotropy for functionally-graded metamaterials \cite{Zhu2021}). Such tools tackle visibility of inner/void regions (mainly in 3D prototypes) to avoid complicated printing paths \cite{Chen2015} and create self-supporting structures \cite{Langelaar2016support} with minimum support material or smooth contouring to avoid notches and acute angles. Level-set \cite{Wang2003,vanDijk2013} and spline-based methods \cite{Seo2010,Qian2013} suit those very purposes.

Remeshing, subdivision and fitting are widespread post-processing steps \cite{Sigmund2007filter}. These tools are often computationally costly and jeopardize the solution's optimality \cite{Lazarov2016}. For that reason, many topology optimization methods embed manufacturing constraints as a part of the iterative process \cite{Langelaar2016} within a constrained optimization framework \cite{Kim2023}. Despite some drawbacks regarding reduced design space and generality \cite{Subedi2020}, the possibilities are plentiful, e.g., auxetic metamaterial design \cite{Andreassen2014}. Interestingly, if volume fraction is left unconstrained, other design parameters can come into play \cite{Bruggi2012,Zhu2020,Dunning2013,ArredondoSoto2021,Koppen2022}. 

\subsubsection{Current interpolation schemes}

Traditional TO methods yield a binary material (black, $\rho = 1$) and void (white, $\rho = 0$) distribution. Intermediate states (gray, $0 < \rho < 1$) frequently need finer profiling strategies \cite{Diaz1995} to define actual limits for material distribution. Relaxation for intervals through power laws is a common solution in stress-based \cite{Yang1996,Le2009} and stress-constrained TO \cite{Holmberg2013,Bruggi2012}. However, it presents several problems, such as stress singularity (addressed by \cite{Bruggi2008}) and the existence of solutions depending upon the parametrization of such power laws. The latter issue can be solved via linear slope approximation \cite{Petersson1998} and/or intervals for sensitivity filters \cite{Kim2023}.

The use of confidence intervals calls for interpolation techniques. However, this cannot be done to material layouts directly (topologies), rather to their associated descriptive matrices (density or stiffness, for instance) as an array of scalars - avoiding singularities and mesh-dependence by non-local methods \cite{Kang2011}. Different densities can also be interpreted as various materials with their particular properties \cite{Bendse1999,Yi2023}, allowing for a more accurate design by getting rid of artificial density penalization. If volume is left unconstrained, material variability is maximized \cite{SaucedoMora2023}.

\subsection{The Carter-Hayes theory}
The idea signing different materials to each density value is coherent with Ashby's law for material selection \cite{Ashby2005} - Equation \ref{eq:carter-hayes}) relating initial $E_0$ and current $E$ Young's moduli through their associated initial $\rho_0$ and current $\rho$ densities:

\begin{equation}
\frac{E}{E_0}=\left(\frac{\rho}{\rho_0}\right)^\gamma
\label{eq:carter-hayes}
\end{equation}

where $\gamma$ is an adiabatic proportionality parameter.

This empirical expression is coherent several studies on the human bone's mechanical behavior \cite{Carter1976,Carter1977,Schaffler1990}. Although common in Engineering applications, the nominal Young's modulus $E_0$ and density $\rho_0$ are mere homogenized assumptions, not to be found in an actual material cell. For that reason, effective properties with physical meaning (like $x_{Phys}$ in top88.m \cite{Andreassen2010}) are preferred.


\subsection{The mechanostat theory for density evolution}
Wolff-Frost's mechanostat theory \cite{Wolff1892,Frost1987} - translated into mechanics by Huiskes \cite{Huiskes1987,Huiskes2000} and Weinans \cite{Weinans1992} -  dictates a linear relationship betwixt mass growth over time (expressed via density $\rho = m/v_f$ for a fixed volume fraction) and strain energy $\Psi = Ku^2/2$ through $B$ as a mechanistic growth parameter (scalar), as in:

\begin{equation}
\frac{d\rho}{dt}=B\,\Psi
\label{eq:mechano}
\end{equation}

This experimental relationship describes bone remodelling, the biological process by which bone cells (osteocytes) are created (synthesized by osteoblasts) or destroyed (resorption by osteoclasts) according to their load-bearing needs (mechanical imperative) \cite{Hughes2010}. This explains bone tissue adaptation to injuries (fracture, implants), age (osteoporosis) and illnesses (osteogenesis imperfecta). Such rationale will be introduced as a computationally cheap embedded filtering strategy with no need for post-processing steps.



\section{Formulation of the effective density correction}

Analougously to the aforementioned mechanostat, the evolution of density and compliance will be intertwined in this proposal. For every target volume fraction, a different starting point applies, giving way to separated interpolation curves converging to a single point where the structure is entirely solid ($\rho=1$). This is not allowed by equation \ref{eq:carter-hayes}, restricted to ($\rho=0, E=0$) as a start. Thus, an effective density $\rho_{ef}$ is defined as the logarithmic update of the real value $\rho_{re}$ given by the mechanostat equation \ref{eq:mechano}:

\begin{equation}
Ln[\rho_{ef}]=Ln[\rho_{re}]+Ln[d\rho_{re}]
\end{equation}

Which can be rewritten as:
\begin{equation}
Ln[\rho_{ef}]=Ln[\rho_{re}\,d\rho_{re}]
\label{rho1}
\end{equation}

Replacing the mechanostat equation \ref{eq:mechano} into equation \ref{rho1}:

\begin{equation}
\rho_{ef}=\rho_{re}\,B\,\Psi\,dt
\end{equation}

Assuming the $\Psi$ independent from density on the grounds of step discretization, it holds:

\begin{equation}
\frac{d\rho_{ef}}{d\rho_{re}}=\alpha\,\Psi
\label{relationship}
\end{equation}

With $\alpha=B\,dt$\\

Using the Carter-Hayes equation \ref{eq:carter-hayes} again:

\begin{equation}
\rho=\rho_0 \, \left(\frac{E}{E_0}\right)^{(1/\gamma)}
\end{equation}

Which, derived with respect to stiffness, gives:

\begin{equation}
\frac{d\rho}{dE}=\frac{\rho_0 \, \gamma}{E_0} \, \left(\frac{E}{E_0}\right)^{(1/\gamma-1)}
\end{equation}

For the particular case of $\gamma=1$, as a simpler approach:

\begin{equation}
\frac{d\rho}{\rho_0}=\frac{dE}{E_0} 
\label{eq:deri}
\end{equation}

For each initial volume fraction, the curves' respective $\rho_0$ and $E_0$ will be denoted as $\rho$ and $E$ from here onwards. The algorithm is now intended to optimize the strain energy density $\Psi$ through variation of the local density $\rho$ related to the Young's moduli $E$. 

This inverse relationship is given by a constant $\zeta$ as $\Psi=\zeta \, /E$. Hence, the derivative of the Young's modulus with respect to the strain energy density yields:

\begin{equation}
\frac{dE}{d\Psi}=-\frac{\zeta}{\Psi^2}=-\frac{E}{\Psi}
\label{eq:Eenergia}
\end{equation}
 
The next step is to substitute the equation \ref{eq:Eenergia} into equation \ref{eq:deri} to obtain:

\begin{equation}
\frac{d\rho}{\rho}=-\frac{d\Psi}{\Psi} 
\end{equation}

And using the $d\rho_{ef}$ instead of the real $d\rho$ (equation \ref{relationship}):

\begin{equation}
\alpha \, \frac{d\rho}{\rho}=-\frac{d\Psi}{\Psi^2} 
\end{equation}

Through the integration of both sides:

\begin{equation}
    \alpha\, Ln[\rho]+\beta=\frac{1}{\Psi}
    \label{eq:integration_interp}
\end{equation}

where $\beta$ is the sum of the integration constants. Considering the equivalence $\Psi = 2c$ ($c = uKu$ being compliance) and $\rho = m/v$, the latter expression can be reshaped into:

\begin{equation}
 \frac{1}{c} = a \ln\left( v\right) + b
    \label{eq:Logarithmic_densification}
\end{equation}

where $a$ and $b$ are constants. More precisely, $b$ corresponds to the lower limit of compliance $c_{min}$ for a maximum volume fraction (and thus, stiffness), i.e. a domain filled with material and no void: $v_{max}=1$. From the former equation, $b = 1/c_{min}$ can be obtained, in which $c_{min}$ can be computed from simple analogy to the continuum.

Substituting $b$ in Equation \ref{eq:Logarithmic_densification} yields $a$, a case-dependent constant defined by the initial volume fraction $v_0$ (considered of infinite compliance) and the minimum compliance $c_{min}$: $a = -\frac{1}{c_{min}ln(v_0)}$. Substituting again in Equation \ref{eq:Logarithmic_densification}, the resulting iterative scheme for volume fraction update is as follows:
\begin{equation}
 v^{i+1} = \left( \frac{1}{c_{min}}-\frac{1}{c^i} \right)  c_{min}\ln\left( v_0\right)
    \label{eq:iterative_scheme}
\end{equation}

Having defined all constants and limits, Equation \ref{eq:iterative_scheme} describes the evolution curves of a given topology from the minimum to the maximum volume fractions while varying in compliance. Any data point along each of the curves represents an existing topology and thus allows for compliance-volume interpolation across different stages given by data points. These curves and the area below them constitute a constrained design dominion, offering an alternative to other more computationally expensive approaches \cite{Bendse1999,Stolpe2001,Zuo2016}. 

This way, multiple volume fractions can provide different compliances and viceversa, on the designer's demand. If volume fraction is left unconstrained, compliance can be an objective on its own \cite{Bruggi2012,Zhu2020} to generate compliant mechanisms \cite{Dunning2013,ArredondoSoto2021,Koppen2022}. Under this approach, the two most important parameters in SIMP (compliance $c$ as the objective, volume fraction $v$ as a restriction) are directly linked, simplifying the process as one unified expression updates both simultaneously where volume becomes a temporary constraint. This formulation is coherent with some instances of logarithmic growth in nature, such as cell reproduction or virus propagation \cite{Triambak2021}.

\section{The erosion filter}

\begin{figure}[H]
    \centering
    \includegraphics[width=0.32\linewidth]{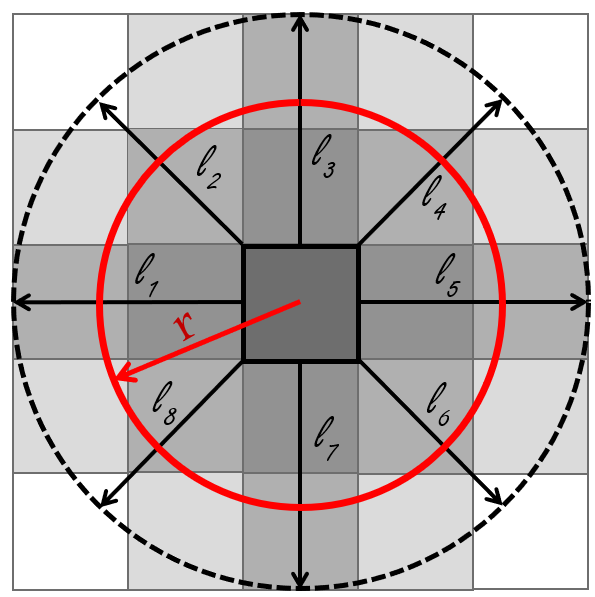}
    \caption{\centering Filter with lengths $l_i, 1 < i < 8 $ (black arrows) outside a radius $r$ (red circle). Each element's color represents its physical density, $\rho_{Phys}$, from white (virtually void) to dark gray (full material).}
    \label{fig:filter_scheme}
\end{figure}

To alleviate the concerns derived from feeble ribs in topology optimization (buckling, stress concentration, fracture, maunfacturability), an additional filter is proposed in this article, besides the sensitivity and density ones already provided by \cite{Andreassen2010}. As an erosion filter, this algorithm searches for neighboring elements with densities below a certain (user-defined) threshold within a given radius $r$ and sets them to void, thus eliminating the unwanted twigs on a minimum-thickness basis. More precisely, it iterates through all non-void elements checking its neighbors at different distances horizontally, vertically and diagonally, in both senses, as in Figure \ref{fig:filter_scheme}. 

Thus, if the advancement in opposite search directions exceeds the target diameter (for instance, $l_1 + l_5 \geq 2r$), the first adjacent elements encountered will be erased if their density is below the established threshold, that is, set to a negligible density value around $0.001$ (non-zero to avoid singularity). Of course, the minimum resolution of this filter is that of the finite element mesh: 1 element. Assuming linear elasticity, scale is left to the user's interest. 

This tool can be considered a minimum-thickness density filter, except its constant effects (applied with the same radius each $X$ iterations) prevent mesh-dependence \cite{Sigmund2007filter}. Other common manufacturability-oriented variables like minimum hole size and sharp edge avoidance are left unrestricted, since the evolution of the filtered topologies is expected to implicitly diminish their effects. The filter can be applied partially, direction-wise, if the user requires so.

After calling the new ad-hoc filter function $twig cutter.m$, the original physical density variable $\rho_{Phys}$ is overwritten with the new filtered value $\rho_{Phys,filtered}$ containing the bi-phasic void/material distribution subject to a user's defined threshold (around $0.8$ for the examples shown in this article). A minimum value is also enforced to avoid checkerboard patterns. This double filter enables local adaptation to manufacturing requirements. 

As an illustrative example, let a 1000x200 cantilever beam be considered, fixed on its left-side wall and loaded with upward unit forces on both its upper and lower right corners. With $top88.m$ \cite{Andreassen2010} as a solver, its embedded density filter is applied with $r_{min} = 2$ (maximum feasible resolution, avoiding checkerboard patterns).

\begin{algorithm}[H]
\caption{Logarithmic density correction and erosion filtering}
\begin{algorithmic}[1]
\Require Setup of the model: boundary conditions and bounding volume.
\Require topological curve limits \( c_{min}, v_{max} \)
\Require logarithmic parameters \( a, b \) with Equation \ref{eq:Logarithmic_densification}
\Require erosion filtering radius \( r \)
\Ensure \textit{FEA isotropic linear elastic calculation at} \( t = 0 \rightarrow K(^{0}E)U = P \) (top88.m \cite{Andreassen2010})
\While{ \( change  > tol \) and \( v  < v_f \)} 
    \State \textit{Initialize} \( \mathbf{c}, \mathbf{dc}, \mathbf{dv}, v_f = v_{min} \)
    \State \textit{Update} \( v \) with Equation \ref{eq:iterative_scheme}
    \For{\( i \) in \( N_{loads} \)}
        \State \textit{Compute} \( K, sK \) (top88.m \cite{Andreassen2010})
        \State \textbf{Initialize} \( (\mathbf{c_i}, \mathbf{dc_i}) = \mathbf{0}\)
        \State \textbf{Initialize} \( \mathbf{dv_i} = \mathbf{1} \)
        \State \textbf{Filter} \( \mathbf{dc_i} \) (sensitivity) and \( \mathbf{dv_i} \) (density)
        \State \textbf{Apply erosion filtering} (twigcutter.m) substituting $\rho_{Phys}$ with $\rho_{Phys,filtered}$
    \EndFor
\EndWhile

\end{algorithmic}
\label{alg:log_filter_code}
\end{algorithm}

\begin{figure}[H]
    \centering
    \includegraphics[width=0.80\linewidth]{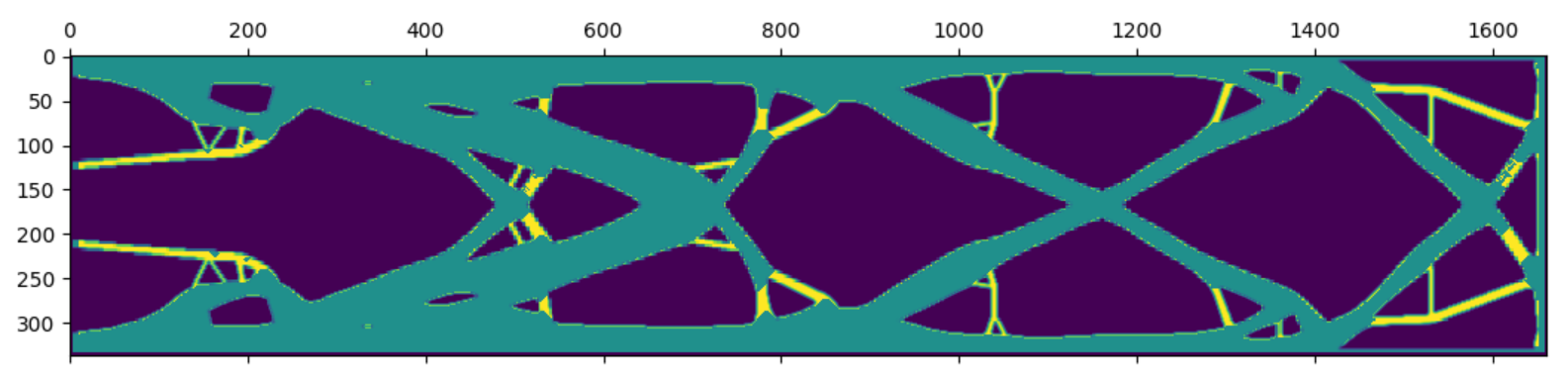}
    \caption{\centering Filter of radius $r=5$ applied to the a 1000x200 cantilever at iteration 1000 (Figure \ref{fig:comparison_it1000}a). Dark blue represent void, yellow the erased (filtered) twigs and light blue the remaining structure.}
    \label{fig:filtering_example}
\end{figure}

Many thin ribs appear in the optimized topology - see Figure \ref{fig:filtering_example}.
This means the distribution of load will be quite heterogeneous: both the normalized shear strain (Figure \ref{fig:comparison_it1000}b) and the von Mises equivalent stress (Figure \ref{fig:comparison_it1000}c) are highly concentrated in junctions, cross-section changes and notches, which puts the structure at risk of fatigue failure (fracture). To get a simpler topology, the filter with radius $r = 5$ is applied in Figure \ref{fig:filtering_example}. 

The applied filter in Figure \ref{fig:filtering_example} effectively tackles most of those structurally superfluous yet difficult to manufacture twigs, which is among the main objectives. Alas, if such a trimming were to be directly introduced, it would not be structurally optimal either: many junctions would be left unsupported, discontinuities would further complicate manufacturing (needing supports if 3D-printed) and the target volume fraction would not be respected after getting rid of those elements. 

Also, some junction points (joints) would be left with undesirable stumps unreachable by the filter as the distance to the main structure is greater than the prescribed radius $r$. Therefore, the filter must be embedded and applied not once at the end (post-processing) but periodically during the iterative optimization to allow for restructuring. 

The filter cannot be enforced during the first iterations as the optimized structure is in its early growth stages, and so densities $\rho_{Phys}$ have not yet reached values close to 1 (full), 0 (void) or any of the upper or lower thresholds to approximate both limits, mostly presenting an ill-defined intermediate density distribution. 

Thus, it would create incompatibilities with the sensitivity/density filters in \cite{Andreassen2010} or even force the full annihilation (void) of the emerging optimized topology, since no existing element's density would reach the established threshold. Alas, choosing the starting iterative step could prove difficult as different loading states generate varied optimized topologies, each of them with their own complexity and time (iterations) to converge.

A side effect of the filtering strategy is that of the decreased volume fraction, which could be leveraged as savings in material. Since the trimming operation takes an unforeseen amount of volume each time, the algorithm needs time to replenish the lost fraction to attain the target. However, if the iteration interval between prunings is too short, the topology might not have the chance to reach the desired volume fraction until no more twigs under the filtering radius are forming, which could take long. Thus, the timing of the filter is also a tuneable parameter. The interval (each $X$ iterations) must be small enough to avoid unwanted regrowth but at the same time sufficiently big so as not to thwart the natural course of topology evolution.

Other option would be leaving volume fraction unrestrained along the iterative process, as proposed in the previous section. To reconstruct a more robust topology after the filter is applied, volume fraction must grow quasi-monotonically - as in equation \ref{eq:iterative_scheme} -, compensating the filter's momentary descent. Otherwise, the filtered structure would be even further restrained and the optimization process would be halted as the only possible change would be continuous retraction of material till total annihilation (full void domain). The practical implementation of the logarithmic densification (Equation \ref{eq:Logarithmic_densification}) and the proposed erosion filter can be seen in Algorithm \ref{alg:log_filter_code}.


\section{New optimization strategy using the growth space}
We propose to distinguish between 3 design spaces for topology optimization designs, depending on the input data and expected outcomes, characterized by four main variables described $\eta$ (set of boundary conditions), $\rho$ (density of the design), $\alpha$ (final topology of the solution), and $\Psi$ (strain energy).

\begin{figure}[H]
    \centering
    \includegraphics[width=0.90\linewidth]{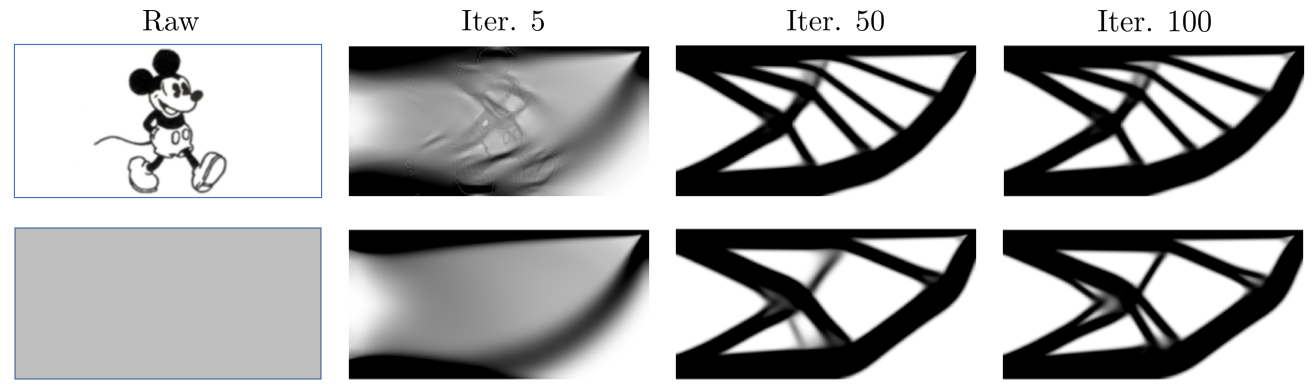}
    \caption{\centering Topology optimization processes (rows) for a cantilever beam with different starting points: a randomly convoluted topology $\alpha$ (top) and the usual intermediate density block ($\rho$ = 0.5, bottom). The same iterations (column-wise) share their boundary conditions $\eta$, volume fraction $\rho$ and strain energy $\Psi$. }
    \label{fig:Mickey}
\end{figure}

With this we can define three spaces according to the optimized designs in the literature. Traditional TO is done in the $\Phi\left(\eta,\rho\right)$ space, with density and boundary conditions as inputs. In this paper, we propose the design space $\Gamma\left(\eta,\rho,\Psi,\alpha\right)$, where we can also impose the desired strain energy of the final structure. We see in figure \ref{fig:Mickey} that a point in the $\Gamma$ space can be reached through different topologies. It gives us insight about a richer space $\Lambda\left(\eta,\rho,\Psi\right)$ that also depends on the topology $\alpha$, which can be prescribed as well. In this work, we present results on the $\Gamma$ space.

This distinction between design spaces is relevant, since any point identified in $\Phi$ by its boundary conditions $\eta$ and volume (density) fraction $\rho$ could be associated with many strain energy states $\Psi$ in $\Lambda$, which in turn correspond to a myriad of possible topologies $\alpha$ in $\Gamma$. This surjection is exemplified in Figure \ref{fig:Mickey}: in the $\Phi$ and $\Lambda$ spaces the top and bottom results are exactly the same point, since both have the same volume fraction and strain energy. The only variation occurs in the $\Gamma$ space due to the different topology $\alpha$ of both results.

\begin{figure}[H]
    \centering
    \includegraphics[width=0.9\linewidth]{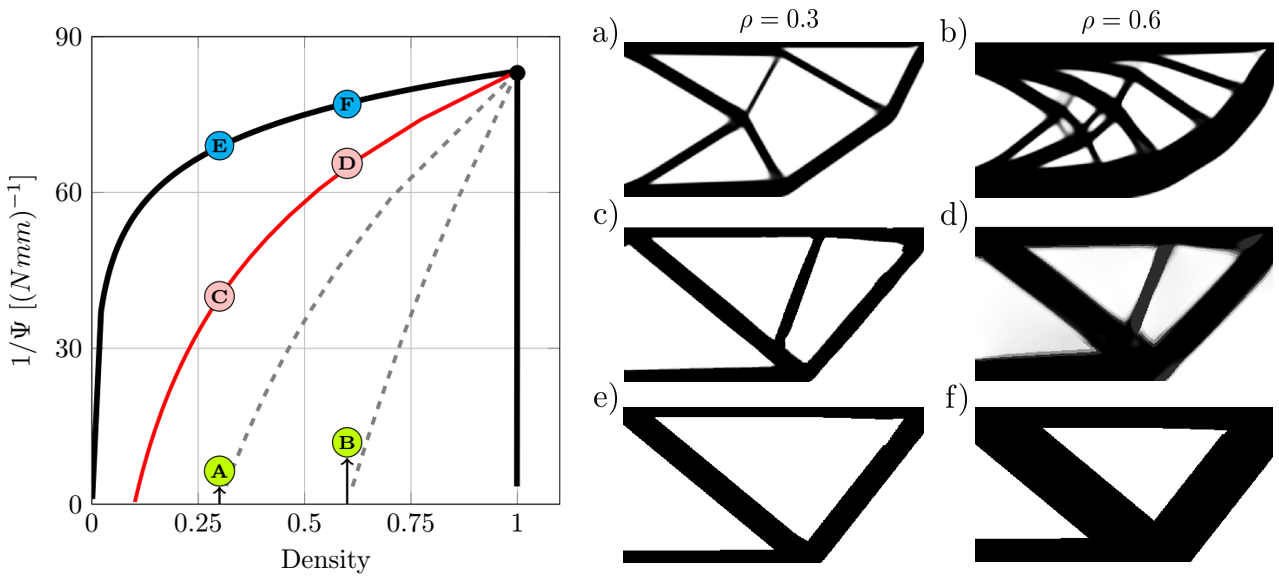}
    \caption{\centering Different optimized cantilever beams (right) with column-wise equal densities $\rho$ represented in the $\Lambda$ space (left). Comparison between SIMP (green dots at the end of vertical arrows) and the proposed logarithmic densification starting from quasi-void (blue dots on the left-side black thick curve) and a very low density ($\rho_0 = 0.1$, pink dots on the red curve).}
    \label{fig:traditional_vs_log}
\end{figure}

However, the starting point for the optimization process is not the only nuance to be considered. Through the proposed logarithmic densification (Equation \ref{eq:Logarithmic_densification}), it is possible to reach equivalent points in $\Phi (\eta, \rho)$ yielding very different strain energies $\Psi$ in $\Lambda$ stemming from distinct starting points. See Figure \ref{fig:traditional_vs_log} for an illustration: since each logarithmic curve designates a different iterative evolution, distinct topologies are expected for the same vertical (density) line as well (designs $A$-$C$-$E$ for $\rho = 0.3$ and $B$-$D$-$F$ for $\rho = 0.6$). 

\begin{figure}[H]
    \centering
    \includegraphics[width=0.9\linewidth]{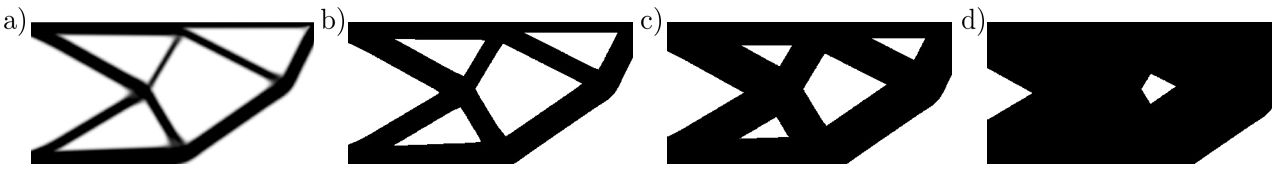}
    \caption{\centering Topological dilation on the $\rho_0 = 0.3$ curve from Figure \ref{fig:traditional_vs_log}. Other examples can be found in Figures \ref{fig:3PB_loginterp} and \ref{fig:f2-r1_dataset_PINN}.}
    \label{fig:dilate}
\end{figure}

Whereas points $A$ and $B$ in Figure \ref{fig:traditional_vs_log} are the result of regular fixed-volume compliance minimization (SIMP: convergence is given by a tolerance between iterations), points $C$-$D$ and $D$-$F$ are obtained via Equation \ref{eq:Logarithmic_densification}, effectively dilating \cite{Sigmund2007} an initial canonical topology: $D$ and $F$ are just thicker versions of $C$ and $E$, respectively. For instance, following the curve corresponding to the initial density for topology $A$, its logarithmic evolution is displayed in Figure \ref{fig:dilate} by dilation of said starting point ($a$). This examples demonstrate $\Lambda$ space's greater versatility when compared to $\Phi$.


\begin{figure}[H]
    \centering
    \includegraphics[width=0.8\linewidth]{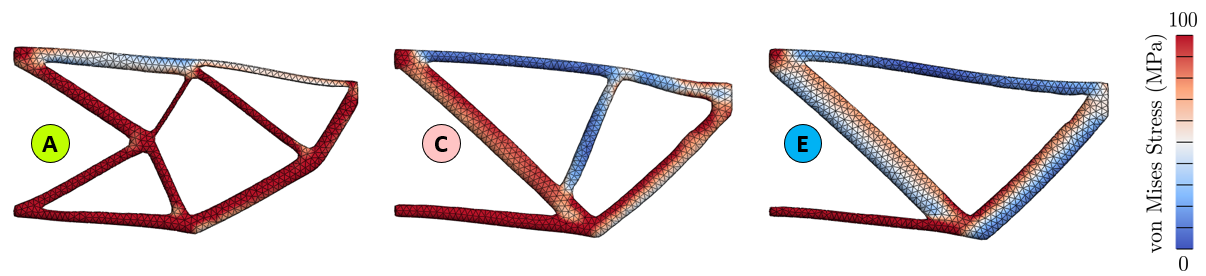}
    \caption{\centering Von Mises stress comparison for topologies $A$, $C$ and $E$ in Figure \ref{fig:traditional_vs_log}.}
    \label{fig:febio_alt}
\end{figure}

The variety of material layout choices available in $\Lambda$ for the very same volume fraction (unobtainable through plain SIMP) has mechanical implications as well. Consider the von Mises stress fields for $A$, $B$ and $C$ (from Figure \ref{fig:traditional_vs_log}) in Figure \ref{fig:febio_alt}: the topological simplification in the proposed designs $C$ and $E$ simpler (w.r.t. $A$, in increasing order) is due to the strengthening and/or disappearance of thin ribs responsible for most of the topology's deformation (and so, compliance) - which can be further reinforced through the proposed erosion filter. This makes the structure stiffer, demonstrated by the lower stress concentration.

\begin{figure}[H]
    \centering
    \includegraphics[width=0.9\linewidth]{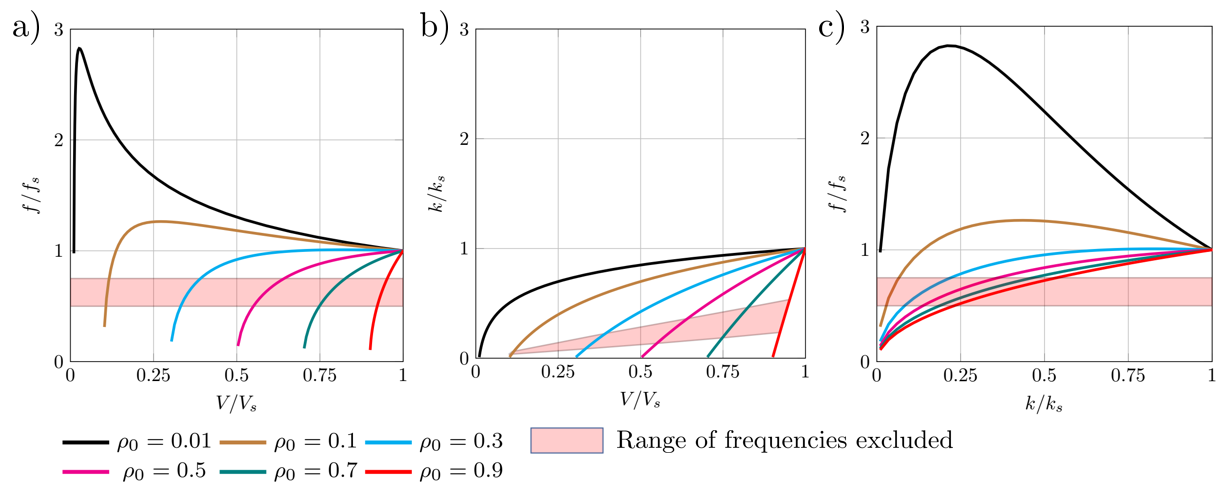}
    \caption{\centering Normalized design curves (w.r.t. solid block values) for volume-frequency (a), volume-stiffness (b) and stiffness-frequency (c) for several initial densities $\rho_0$, showcasing a band-gap in red.}
    \label{fig:frequency_vs_k_v}
\end{figure}

\section{Design for a target natural frequency}

For a sole isotropic material, i.e. SIMP's assumptions, parameters $\rho$ and $\Psi$ within the $\Lambda$ design space are tantamount to mass $m$ and stiffness $k$, respectively - allowing to compute the prototype's natural frequency $f_0$:

\begin{equation}
    f_0 = \frac{1}{2\pi}\sqrt{\frac{k}{m}}
    \label{eq:natural_frequency}
\end{equation}

This offers many direct applications in Engineering, such as vibration prevention or acoustic resonator design \cite{Wang2017_res} - a common goal for metamaterials. Although SIMP is by default designed to maximize stiffness, the regular approach \cite{Bendse1988} does not allow for volume (or density) tailoring and yields artificial frequency modes when low-density cells are penalized. Hence the use of genetic/evolutionary methods by some authors \cite{Picelli2015,Huang2010}, where volume is unconstrained. 

The iterative scheme suggested in this article (Equation \ref{eq:iterative_scheme}) permits flexible volumes and frequency control without solving costly eigenvalue problems \cite{Leader2019}. For instance, to select (or avoid) a frequency value or range, it suffices to visualize the frequency/stiffness/volume interval and its intersections with predefined topological curves and select iterative parameters $\rho_0$, $\rho_i$, $v_i$ accordingly. See Figure \ref{fig:frequency_vs_k_v} for an example.

\section{Results}

In this section, the proposed tools (volume-compliance interpolation, erosion filtering and their respective and combined applications) will be put to the test on several practical case studies.

\subsection{Filtering radius in topology optimization}

The filtering algorithm $twigcutter.m$ will be called every 10 iterations, considering it a lapse long enough to polish the protruding remains while short enough to impede any pernicious regrowth of the freshly cut parts. To avoid further constriction of the design space, the volume fraction grows constantly by a fixed amount starting from an initial prescribed value (about 90\% of the target) until the target fraction is met. 

This way, the structure has room to regrow structurally meaningful ribs while obtaining the desired target as the twigs are periodically severed. Sometimes, the eliminated volume is greater than what can be regenerated within the filtering interval, so the volume fraction is somewhat lower than the desired target - which could be construed as an additional advantage, since material efficiency is further ensured. 


\begin{figure}[H]
    \centering
    \includegraphics[width=1\linewidth]{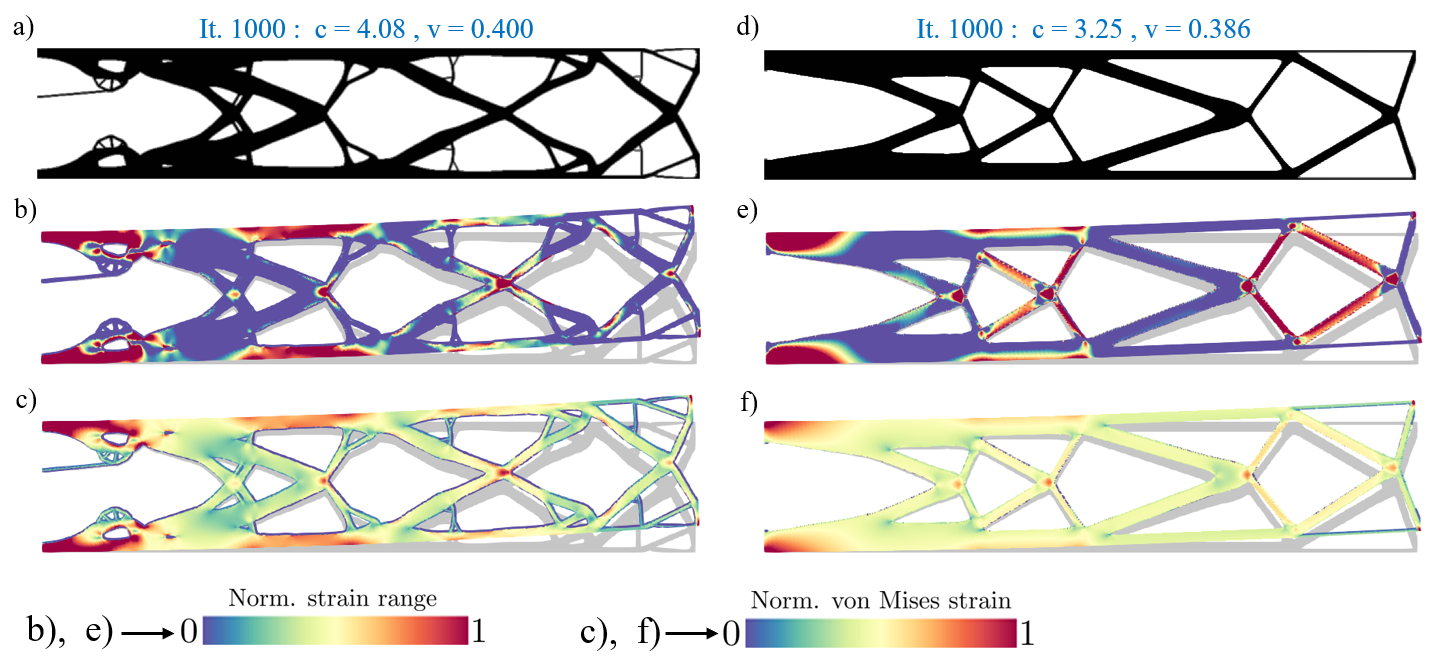}
    \caption{\centering Vanilla (left) and filtered (right) topology optimization of a cantilever beam with upward loads in both right corners, iteration 1000. The resulting optimized structures (a, d) are shwon along their shear strain and von Mises equivalent stress distributions (b,c and e,f respectively).}
    \label{fig:comparison_it1000}
\end{figure}

The topologies obtained after 1000 iterations are showcased in Figure \ref{fig:comparison_it1000}, with their respective normalized shear strain and equivalent von Mises stress ($b$ and $c$ for the unfiltered version ($a$), $e$ and $f$ for the filtered one ($d$)). Focusing on the left side (Figure \ref{fig:comparison_it1000}a, b and c), it is easy to notice how the thinnest struts are not actually performing any strain or stress transfer worth of mention in the vanilla topology - especially regarding shear strains since loading is vertically applied -, rendering them useless. As a result, both shear and strain are highly concentrated in some joints and near loading points and supports, whereas the remainder of the structure is idle, supporting a fraction of the load. 


The filtered version (see Figure \ref{fig:comparison_it1000}d) is noticeably simpler and sturdier: the genus, i.e., “number of holes”, has decreased greatly and the majority of those inconvenient twigs have disappeared giving way to thicker, more robust ones with a meaningful contribution to load bearing - namely between joints and support points. Additionally, the stumps seen in Figure \ref{fig:filtering_example} have been completely erased as a result of the iterative process. 


This has noticeable effects on their shear strain (Figure \ref{fig:comparison_it1000}, second row) and von Mises equivalent stress (Figure \ref{fig:comparison_it1000}, third row): having fewer struts guarantees a more consistent and homogeneous load distribution, so that a smaller set of ribs are structurally meaningful, instead of unevenly loaded sections of struts which create unnecessary and dangerous stress concentrations. This is preferable both from a mechanical and a practical point of view (manufacturability), since reinforcing those individual struts entirely (e.g., making it thicker or choosing a sturdier material) is easier than doing so with multiple local points.

\begin{figure}[H]
    \centering
    \includegraphics[width=1\linewidth]{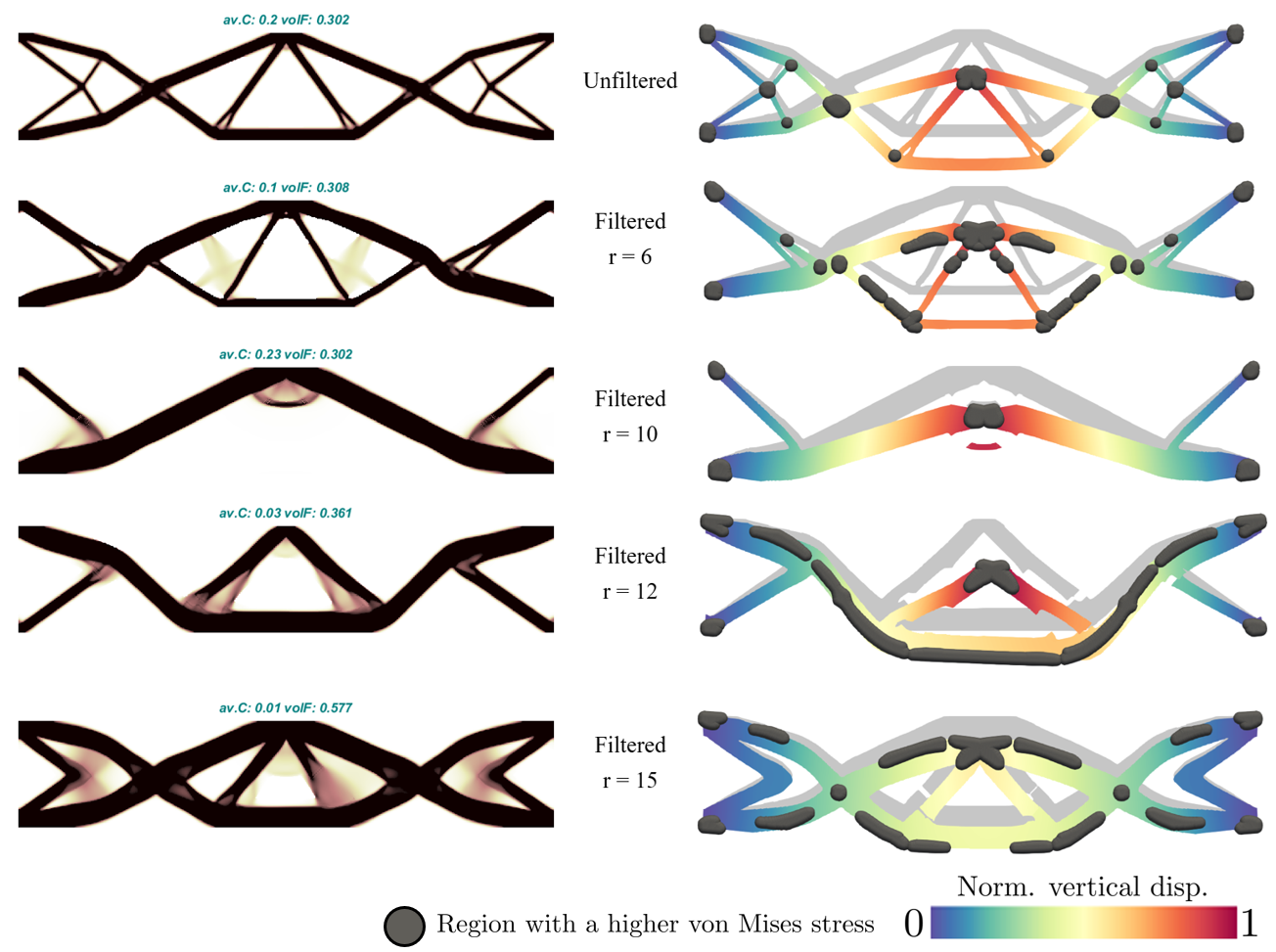}
    \caption{ \centering Optimized topologies for a beam fixed on both ends under a unitary downward load in the middle upper section (left) with their respective vertical displacement (color map) and von Mises (over 0.05 MPa in gray) equivalent stresses on the right: unfiltered ($1^{st}$ row), $r=6$ ($2^{nd}$ row), $r=10$ ($3^{rd}$ row), $r=12$ ($4^{th}$ row), $r=15$ ($5^{th}$ row). }
    \label{fig:bifixed_various_radii}
\end{figure}

Interestingly, zero-shear areas (in blue, second row in Figure \ref{fig:comparison_it1000}) imply equal principal stresses ($\varepsilon_1 = \varepsilon_2$), which could be further leveraged for manufacturing purposes, e.g. by printing filaments along those principal directions so strength is maximal.

On top of the mechanical improvement, computation power is also well administered: since the trimming process is done at the same time as the regular topology optimization, it does not imply any noticeable additional computational power, being run on a single core i7 1.8 GHz, 16 GB RAM in about 1.5h time till iteration 1000.

The filter is suitable for any applied loads and boundary conditions. For instance, a double-fixed beam with a single unitary downward load is subjected to filtered topology optimization with various radii, whose results can be seen in Figure \ref{fig:bifixed_various_radii}. As expected, topologies get simpler on the filtered cases, although not always monotonically with an increasing radius, since boundary conditions still play an important role. 

As the radius grows, two mutually linked effects are noticeable (see Figure \ref{fig:bifixed_various_radii} left, on top of each topology): volume fraction grows and compliance decreases (that is, stiffness increases). The former is expected, considering the algorithm needs to restock structurally important ribs in an ever-thicker version to comply with filtering requirements, which in turn provokes compliance descent: topologies become simpler (fewer deformable parts) and stiffer (less deformation overall). This can be clearly seen in their respective vertical displacement color maps in Figure \ref{fig:bifixed_various_radii} right.

As hinted before, von Mises equivalent stress (highest values in grey in Figure \ref{fig:bifixed_various_radii} right) is homogenized in filtered versions: continuous areas along most structurally relevant struts as opposed to local points in junctions and supports in unfiltered equivalents. Density in filtered topologies is less sharply defined than in its unfiltered counterpart (Figure \ref{fig:bifixed_various_radii} first row), even less so for greater radii. This boils down to the need to recover greater volumes scrapped off by the filter, which takes longer if the trimmed fraction is bigger, so convergence to the underlying optimal topology is slower.

\begin{figure}[H]
    \centering
    \includegraphics[width=1\linewidth]{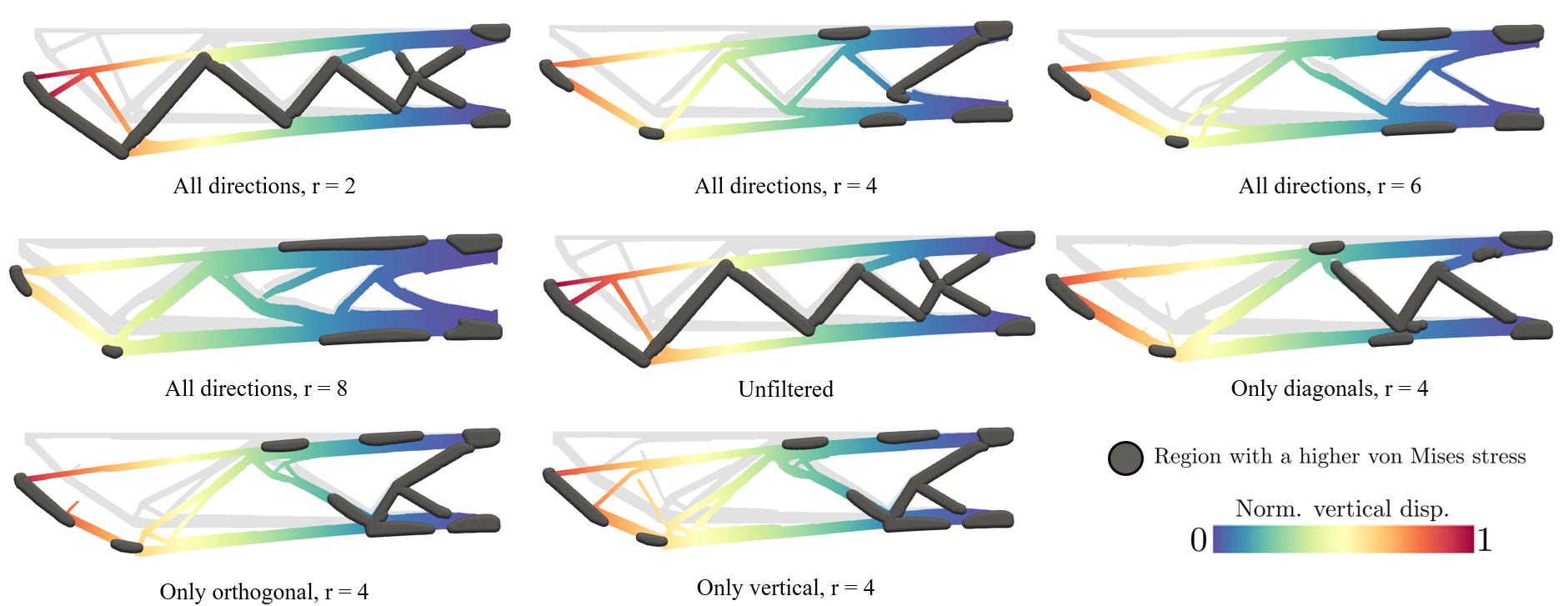}
    \caption{\centering Displacement (color scale) and von Mises stress (over 0.05 MPa in gray) for different optimized cantilever beams (downward unit force on the upper left corner) under erosion filtering.}
    \label{fig:cantilever_multdir}
\end{figure}

The filter can also be applied partially, i.e. exclusively in some search directions. See Figure \ref{fig:cantilever_multdir} for a cantilever example. Again, filtering effects are easier to see when the radius is high: while there's virtually no difference for lower radii (compare $Unfiltered$ and $All, r = 2$, for instance), heavier filtering creates simpler topologies with more homogeneous displacement (color map) and stress (gray areas); overall stiffer (notice the smallest tip displacement for the biggest filtering radius, $r = 8$). 


However, asymmetric filtering can yield different results. Diagonally applied filters favor diagonal struts, longer in the center section and mostly stressed near the support and loading point (see \ref{fig:bifixed_various_radii} center right for $r = 4$ and directions $135^{\circ}$ and $315^{\circ}$, that is, $l_2$ and $l_6$ in Figure \ref{fig:filter_scheme}). Orthogonality (directions $l_1$, $l_3$, $l_5$ and $l_7$ in Figure \ref{fig:filter_scheme}) offers further ramification in the remaining unfiltered directions (diagonals) and more localized stresses. If the filter is just applied vertically (directions $l_3$ and $l_7$ in Figure \ref{fig:filter_scheme}), horizontal capillarity is allowed and thus some ribs are closed. 

This partial filtering could prove useful to induce tailored anisotropy for mechanical reasons (e.g. printing directions), thus controlling the incidence of buckling, compression, traction, etc. to meet the designer's requirements (material, boundary conditions, reliability, manufacturing method).

\subsection{Structural designs within the growth space}

\begin{figure}[H]
    \centering
    \includegraphics[width=1\linewidth]{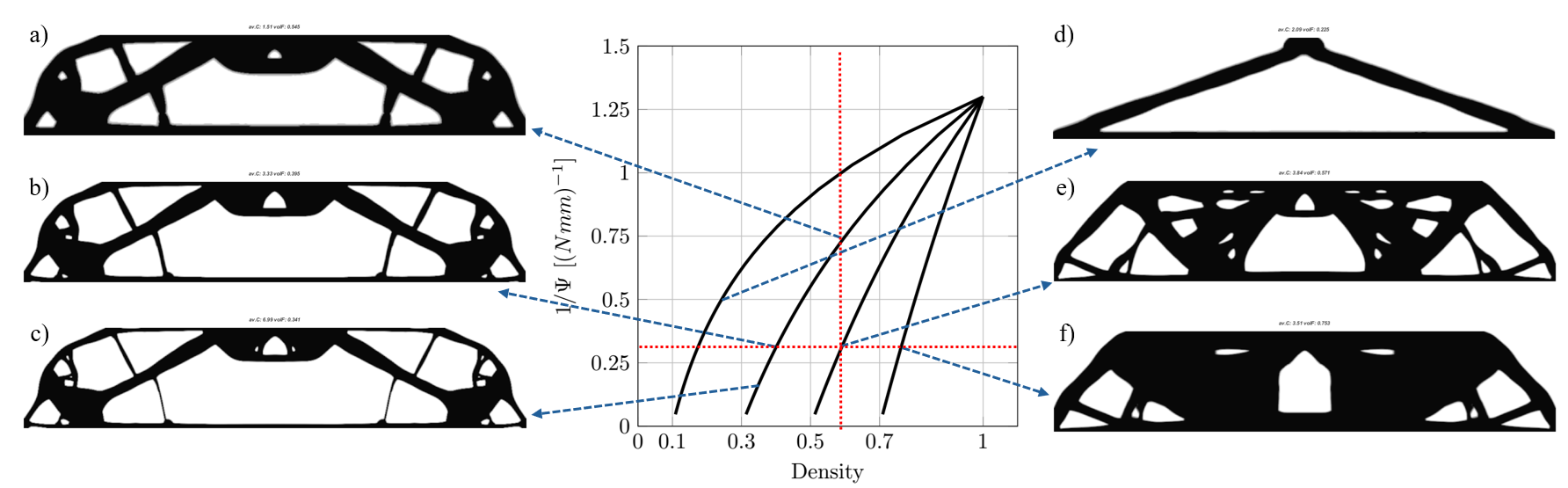}
    \caption{\centering Unfiltered evolving topologies for a 3-point bending structure and several initial volume fractions: $0.1$ (left), $0.3$ (middle left), $0.5$ (middle right) and $0.7$ (right).}
    \label{fig:3PB_nofilter_log}
\end{figure}

An example of logarithmic optimization of a 3-point bending beam with downward unit forces applied in half length with different volume fractions and compliances can be seen in Figure \ref{fig:3PB_nofilter_log}, as a result of several independent optimization processes starting on different minimum volumes $v_0 = 0.1$, $v_0 = 0.3$, $v_0 = 0.5$ and $v_0 = 0.7$. The detailed information for data points $a$-$f$ is contained in Table \ref{tab:Nofilter_log_data}.

\begin{table}[H]
    \centering
    \begin{tabular}{|c|c|c|c|c|c|c|}
        \hline
         \textbf{Fig. 7} & \textbf{a} & \textbf{b} & \textbf{c} & \textbf{d} & \textbf{e} & \textbf{f} \\
        \hline
        \textbf{$v_0$} & 0.3 & 0.3 & 0.3 & 0.1 & 0.5 & 0.7 \\
        \hline
        \textbf{$v_i$} & 0.545 & 0.395 & 0.341 & 0.225 & 0.571 & 0.753 \\
        \hline
        \textbf{$c_i$} & 1.51 & 3.33 & 6.99 & 2.09 & 3.84 & 3.51 \\
        \hline
        $i$ & 65 & 45 & 25 & 125 & 20 & 15 \\
        \hline
    \end{tabular}
    \caption{\centering Information associated to topologies shown in Figure \ref{fig:3PB_nofilter_log}: initial $v_0$, and current volume fraction $v_i$ and compliance $c_i$ for iteration $i$ (3-point bending).}
    \label{tab:Nofilter_log_data}
\end{table}

Focusing on Figure \ref{fig:3PB_nofilter_log}'s left side, the development of a topological family with $v_0 = 0.3$ can be seen at different iterations: 25 (c), 45 (b) and 65 (a). Since the filter has not been applied, the only possible way to grow in volume is widening the already existing ribs, which in turn lowers compliance (structures get stiffer). This approach resembles plain shape optimization, where topology remains practically invariant (aside from little holes being engulfed by strut thickening alone). 

If the designer is aiming at lower volume fractions, the optimized topologies get much simpler (triangular, linking loading points and supports directly) and relatively stiff, although they take longer to converge since the initial volume is too low to evolve into low-compliance structures too soon. See Figure \ref{fig:3PB_nofilter_log}d (starting volume $v_0 = 0.1$) for an example with $v_i = 0.225$ and $c_i = 2.09$ at iteration $i = 125$ (Table \ref{tab:Nofilter_log_data}). 

Should the modeler want stiffer options, starting with higher volume fractions could help speed up convergence to a desired compliance. Consider labels $b$, $e$ and $f$ in Figure \ref{fig:3PB_nofilter_log} and their respective data in Table \ref{tab:Nofilter_log_data}: while all three topologies have approximately the same compliance $c_i$, their respective volume fractions $v_i$ are very different ($f$'s being almost twice as big as $b$'s), and their time to convergence is too ($b$ takes three times as many iterations $i$ to generate compared to $f$).

Thus, $v_0$ becomes a key user-defined hyper-parameter controlling both computation time (number of iterations) and final volume fraction $v_i$ for a given compliance $c_i$. The needs of each topology will be determined by the material and boundary conditions. Interestingly, isochoric (vertical) lines, like the one (approximately) linking topologies $a$ and $e$, represent processes with variation of compliance for the same volume. 

Note that moving upwards in Figure \ref{fig:3PB_nofilter_log} means lowering compliance while maintaining volume fraction, i.e. the regular SIMP method. Going downwards is also possible if the goal is a more compliant structure, among many other strategies neglected by SIMP.

\begin{figure}[H]
    \centering
    \includegraphics[width=1\linewidth]{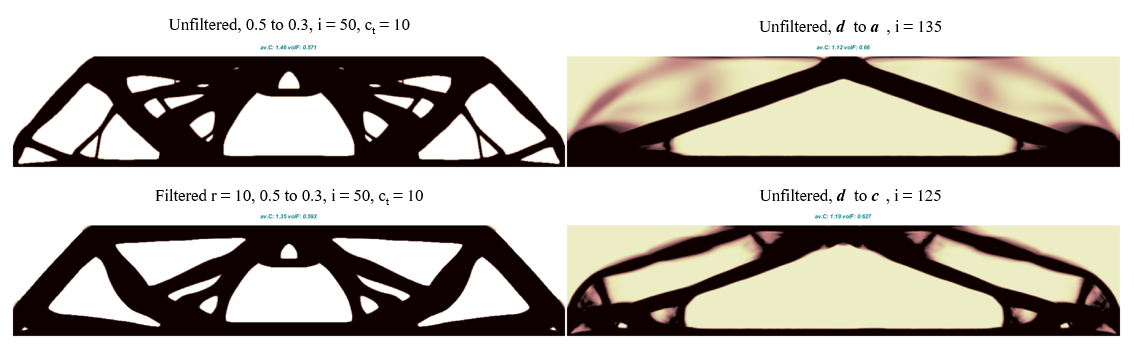}
    \caption{\centering Hybrid topologies generated by interpolation between different $v_0$ curves: mixed 0.3-0.5 topology at iteration 50 with (lower left) and without $r=10$ filter (upper left) with a compliance threshold of $c_t = 10$ and $v_0 = 0.3$ curves and volume-threshold interpolation between $d$ and $a$ (upper right) and between $d$ and $c$ (lower right).}
    \label{fig:3PB_loginterp}
\end{figure}

This way, design interpolation would be possible by jumping between different logarithmic curves (with their respective $v_0$). Transition between curves can be done horizontally (fixing a compliance limit $c_t$ and varying the volume fraction $v_i$ till intersection with the desired volume curve) or vertically (fixing a volume threshold $v_t$ and adjusting the compliance $c_i$ until the objective curve is reached). See Figure \ref{fig:3PB_loginterp} for some examples of interpolation (filtered/unfiltered, horizontal/vertical).

On Figure \ref{fig:3PB_loginterp}'s left side, some examples of horizontal interpolation (compliance threshold $c_t = 10$) are displayed. Both topologies closely resemble the initial curve's family ($v_0 = 0.5$, e.g. Figure \ref{fig:3PB_nofilter_log}e). For the same iteration ($i = 50$), the filtered case ($r = 10$, lower left) yields a simpler material layout than the unfiltered case (upper left), a “cleaner” structure with equivalent mechanical properties: similar compliance (up 1.46 vs down 1.35) and volumes (up 0.571 vs down 0.593). An exact match could be found if needed. \\

Figure \ref{fig:3PB_loginterp} right (vertical interpolation with volume threshold $v_t$) showcases the importance of both the initial curve (determined by its initial volume fraction $v_0$) and the amount of iterations on each curve: while the resulting topologies do reflect aspects of both curves, the initial one ($v_0 = 0.1$) clearly prevails, judging by the absence of central vertical ribs (present in both $a$ and $c$) under the load-support triangle (coming from $d$). 

Thus, path-dependence applies to this kind of interpolation schemes, and so it must be taken into account accordingly. These intermediate solutions remain close in compliance (up 1.12 vs down 1.19) and somewhat in volume (up 0.660 vs down 0.627). Reaching the exact targets (Figures \ref{fig:3PB_nofilter_log}a and \ref{fig:3PB_nofilter_log}c, respectively) becomes difficult due to the many factors involved (initial volume, interpolation, threshold, boundary conditions, filters, etc.). Importantly, volume fraction (and thus, density) can only be fixed or increased ($v_{i+1}\geq v_i$), since a decrease in volume contradicts Equation \ref{eq:iterative_scheme}. 

The obvious differences between horizontal ($c_t$) and vertical ($v_t$) interpolation are due to the chosen iterative scheme: according to Equation \ref{eq:iterative_scheme}, the next iteration's volume fraction $v^{i+1}$ is directly depending on the current iteration's (averaged) compliance $c^i$ - an exclusive result of the undergoing SIMP optimization -. This implies that constant compliance interpolation (horizontal, Figure \ref{fig:3PB_loginterp} left) does not incur in any meaningful volume changes (unless a filter is applied, see Figure \ref{fig:3PB_loginterp} lower left). Hence, topologies remain more or less unchanged when compared to their initial curve's ($v_0 = 0.5$) canonical shape (same “topological family” as Figure \ref{fig:3PB_nofilter_log}e).\\

\begin{figure}[H]
    \centering
    \includegraphics[width=1\linewidth]{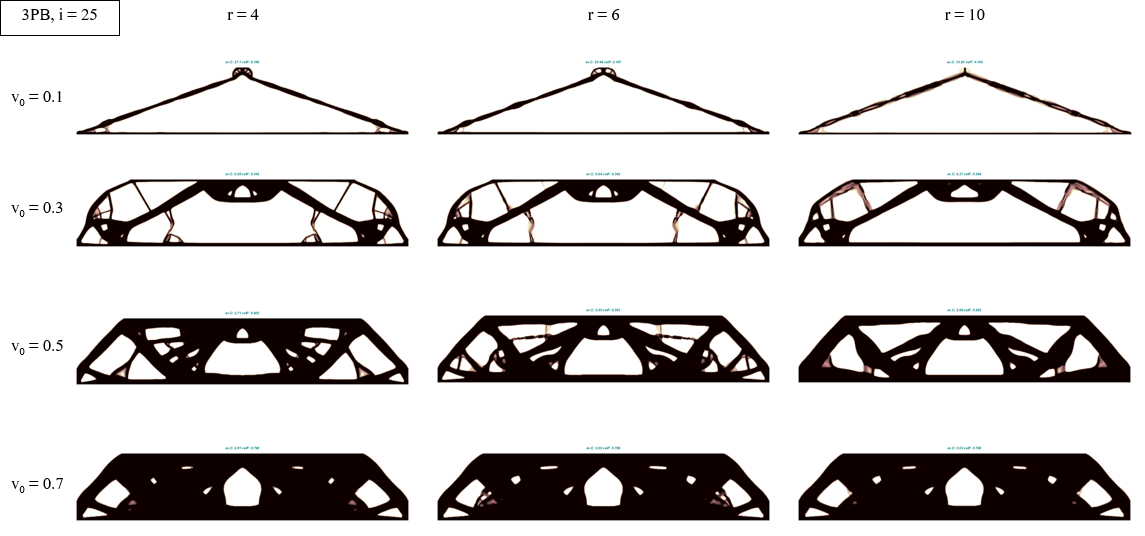}
    \caption{\centering Filtered versions of the 3-point bending topology with different radii (4, 6 and 10) and starting volume fractions (0.1, 0.3, 0.5 and 0.7). Iteration 25. }
    \label{fig:3PB_filter_log_it25}
\end{figure}

\begin{table}[H]
    \centering
    \begin{tabular}{|c|c|c|c|}
        \hline
         \textbf{$c_i | v_i$ (i = 25)} & \textbf{r = 4} & \textbf{r = 6} & \textbf{r = 10} \\
         \hline
         \textbf{$v_0 = 0.1$} & 27.1 | 0.106 & 25.44 | 0.107 & 33.83 | 0.104 \\
         \hline
         \textbf{$v_0 = 0.3$} & 6.58 | 0.342 & 6.64 | 3.42 & 6.37 | 0.344 \\
         \hline
         \textbf{$v_0 = 0.5$} & 2.71 | 0.603 & 2.95 | 0.593 & 2.98 | 0.593 \\
         \hline
         \textbf{$v_0 = 0.7$} & 2.01 | 0.798 & 2.02 | 0.798 & 2.02 | 0.798 \\
         \hline
    \end{tabular}
    \caption{\centering Compliance $c_i$ and volume fraction $v_i$ for topologies shown in Figure \ref{fig:3PB_filter_log_it25}.}
    \label{tab:Tab_volcomp_3PBit25}
\end{table}

Conversely, vertical interpolation (constant volume $v_t$, Figure \ref{fig:3PB_loginterp} right) allows a free evolution of compliance between curves, so an intermediate trade-off solution is found between the beginning and ending “canonical topologies”, with a greater influence of the latter for already disclosed reasons. In these cases, convergence is slow and ill-defined.  Interpolation between more than two curves would also be possible, although complex, needing enough iterations in between.

Leveraging both the filter and the logarithmic approach, very diverse topologies can be obtained for the very same design requirements (compliance, volume fraction). Let the 3-point bending beam be considered again. Several starting volumes $v_0$ and filtering radii $r$ are enforced to produce a plethora of design points with similar compliance $c_i$ and volume fraction $v_i$ but very different topologies. 

Figure \ref{fig:3PB_filter_log_it25} and Table \ref{tab:Tab_volcomp_3PBit25} show an array of topologies for iteration $i = 25$ and various starting volumes $v_0$ (0.1, 0.3, 0.5, 0.7) and filtering radii $r$ (4, 6, 10), with their respective compliances $c_i$ and volume fractions $v_i$. 

Row-wise, greater filtering radii (rightwards) “simplify” the topologies, reducing their $genus$ (“number of holes”), as observed in previous examples. Column-wise, bigger starting volume fractions $v_0$ (downwards) densify the structures, performing a sort of shape optimization (as in Figure \ref{fig:3PB_nofilter_log}$c$ to $a$) while significantly altering their genus. 

\begin{figure}[H]
    \centering
    \includegraphics[width=1\linewidth]{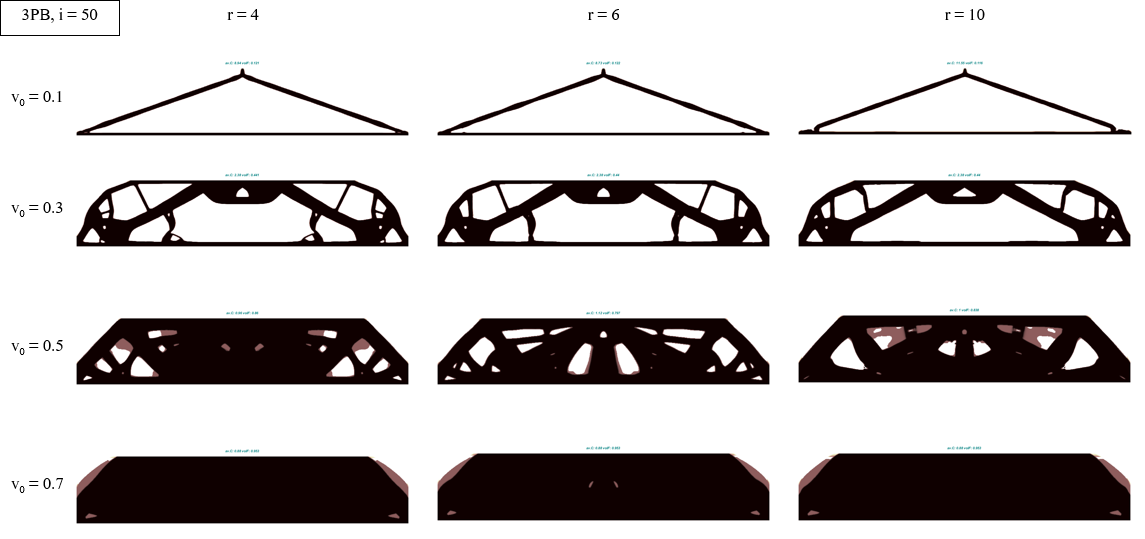}
    \caption{\centering Filtered versions of the 3-point bending topology with different radii (4, 6 and 10) and starting volume fractions (0.1, 0.3, 0.5 and 0.7). Iteration 50. }
    \label{fig:3PB_filter_log_it50}
\end{figure}

\begin{table*}[ht]
    \centering
    \begin{tabular}{|c|c|c|c|}
        \hline
         \textbf{$c_i | v_i$ (i = 50)} & \textbf{r = 4} & \textbf{r = 6} & \textbf{r = 10} \\
         \hline
         \textbf{$v_0 = 0.1$} & 8.94 | 0.121 & 8.73 | 0.122 & 11.55 | 0.116 \\
         \hline
         \textbf{$v_0 = 0.3$} & 2.38 | 0.441 & 2.38 | 0.440 & 2.38 | 0.440 \\
         \hline
         \textbf{$v_0 = 0.5$} & 0.96 | 0.86 & 1.12 | 0.797 & 1.00 | 0.838 \\
         \hline
         \textbf{$v_0 = 0.7$} & 0.88 | 0.953 & 0.88 | 0.953 & 0.88 | 0.953 \\
         \hline
    \end{tabular}
    \caption{\centering Compliances $c_i$ and volume fractions $v_i$ for topologies shown in Figure \ref{fig:3PB_filter_log_it50}.}
    \label{tab:Tab_volcomp_3PBit50}
\end{table*}

According to Table \ref{tab:Tab_volcomp_3PBit25}, topologies in the same row (same $v_0$, different $r$) present very similar volume fractions $v_i$. This is expected, since they are contained in the same $v_0$ curve under the same number of iterations $i$. 

Such is the case of their compliances ($c_i$), thus representing different topologies for virtually the same data point ($v_i$, $c_i$) in Figure \ref{fig:3PB_nofilter_log} - which offers an unmatched design flexibility. This row-wise equivalence is truer the closer radii are to each other and the greater the starting volume fractions are, where the design space is more constrained and coincidence is almost exact. 

For very low starting volumes (e.g., $v_0 = 0.1$), their simpler topology is far less limited and so distinct features for each filtering radius provoke relatively important fluctuations (see middle top on the topologies shown in Figure \ref{fig:3PB_filter_log_it25} first row) in volume and, most importantly, in compliance. Nevertheless, these discrepancies can be solved by slightly varying iterations, since it has been mentioned that convergence follows different paces depending on the chosen logarithmic curve and filtering strategy.


Observing Figure \ref{fig:3PB_filter_log_it50} and Table \ref{fig:3PB_filter_log_it50} for iteration 50, some of the topologies seen in Figure \ref{fig:3PB_filter_log_it25} have noticeably evolved (middle volumes, namely $v_0 = 0.3$ and $v_0 = 0.5$), while others remain practically unaltered (lower volume, $v_0 = 0.1$). Higher volumes ($v_0 = 0.7$) have almost reached the maximum possible stiffness (highest volume, lowest compliance), i.e. the full material block. Table \ref{tab:Tab_volcomp_3PBit50} reflects closer compliance and volume values for the same row than Table \ref{tab:Tab_volcomp_3PBit25}, virtually identical for $v_0 = 0.3$ and $v_0 = 0.7$.

All the remarks and observations applied for this case of study remain true for different loads, boundary conditions, radii and starting volumes, as soon as linear elasticity is kept. The combination of both tools presented in this article (filtering and logarithmic densification) provides a powerful and versatile inverse design methodology with three degrees of freedom: compliance $c_i$, volume fraction $v_i$ and filtering radius $r$ (somewhat equivalent to minimum thickness).

\section{D$^2$NN interpolation with effective density correction}

This example is designed to implement a Double Distance Neural Network (D$^2$NN) with Equation \ref{eq:iterative_scheme} as a physical loss function evaluating how far is the evaluated data point from the corresponding logarithmic curve, as well as the distance within. 

\begin{algorithm}[H]
\caption{D$^2$NN interpolation}
\begin{algorithmic}[1]
\Require Setup of the model: boundary conditions $\eta$ and bounding volume $\rho$.
\Require topological curve limits \( c_{min}, v_{max} \)
\Require logarithmic parameters \( a, b \) with Equation \ref{eq:Logarithmic_densification}
\Require erosion filtering radius \( r \)
\Ensure \textit{FEA isotropic linear elastic calculation at} \( t = 0 \rightarrow K(^{0}E)U = P \) (top88.m \cite{Andreassen2010})

\State \textbf{Define} function \textit{resize image} to scale image resolution
\State \textbf{Define} function \textit{load data} for \( c, v \) and topology image files
\State \textbf{Define} class \textit{D$^2$NN} featuring a sequential architecture and forward propagation
\State \textbf{Define} function \textit{D$^2$NN loss} with Equation \ref{eq:iterative_scheme}.
\State \textbf{Define} function \textit{D$^2$NN training} applying the \textit{adam} optimizer for N epochs
\State \textbf{Define} function \textit{generate image} from \( c, v \) targets
\State \textbf{Define} function \textit{save image}

\State \textbf{Call} function \textit{load data} upon the (\( c, v \), image) dataset
\State \textbf{Call} function \textit{D$^2$NN training} to train the \textit{D$^2$NN} model with the loaded dataset
\State \textbf{Set} design targets \( c_t, v_t \) 
\State \textbf{Call} function \textit{generate image} upon the design targets \( c_t, v_t \) 
\State \textbf{Call} function \textit{save image} upon the new output image

\end{algorithmic}
\label{alg:ML}
\end{algorithm}

Therefore, an interpolation framework is created to efficiently define topologically optimized structures with a given volume fraction and strain energy within a design space $\Lambda(\eta,\rho,\Psi)$, bypassing the need for relaxation methods or multi-material combinations anf their inherent computation complexity.

A D$^2$NN-driven surrogate is devised to accelerate TO convergence, featuring alternating linear (3) and ReLU layers (2) and a final sigmoid pooling. The model is trained for 500 epochs to predict compliance and volume values from topologies as image files, as in Algorithm \ref{alg:ML}. A 600x200 cantilever beam is considered for a case study. 

\begin{figure}[H]
    \centering
    \includegraphics[width=1\linewidth]{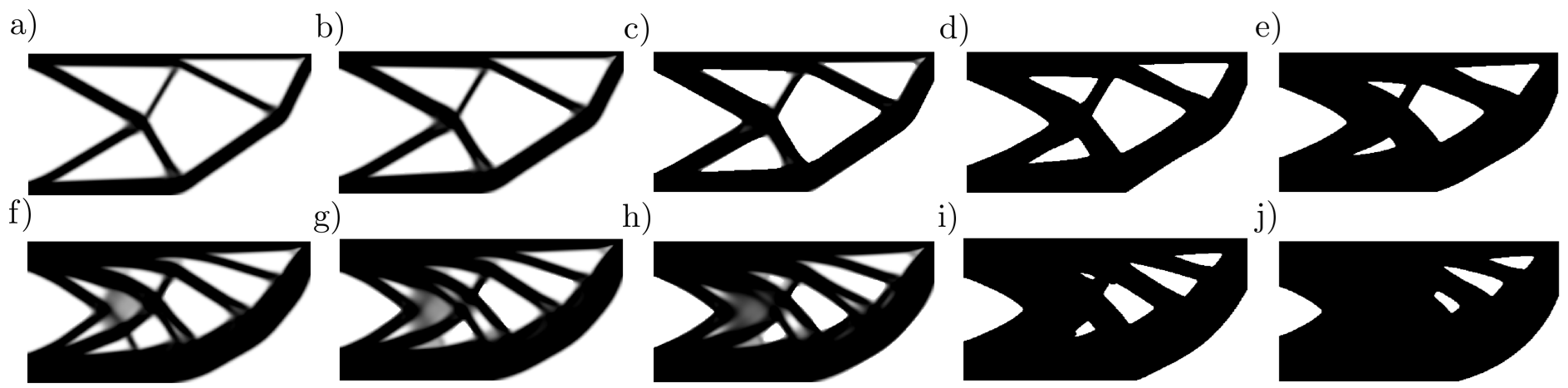}
    \caption{\centering Training set: points A-J in Figure \ref{fig:f1-curves_log} for the ML surrogate with filtering radius $r = 1$.}
    \label{fig:f2-r1_dataset_PINN}
\end{figure}

\begin{figure}[H]
    \centering
    \includegraphics[width=1\linewidth]{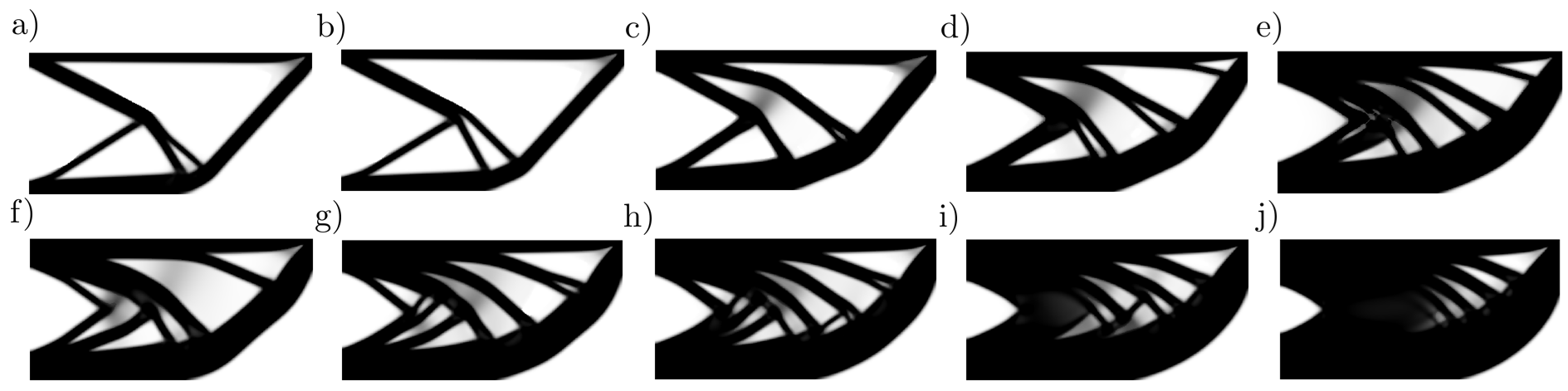}
    \caption{\centering Training set: points A-J in Figure \ref{fig:f1-curves_log} for the ML surrogate with filtering radius $r = 4$.}
    \label{fig:f3-r4_dataset_PINN}
\end{figure}

\begin{figure}[H]
    \centering
    \includegraphics[width=0.45\linewidth]{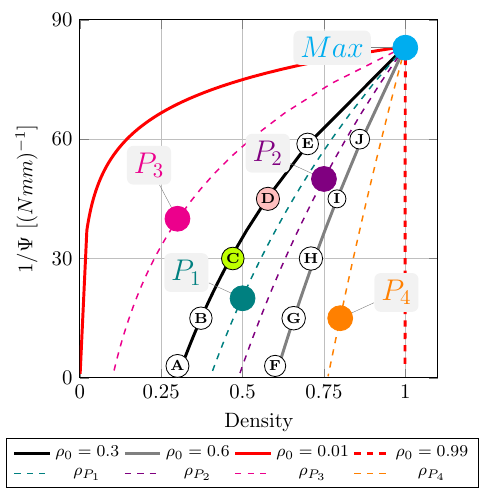}
    \caption{\centering Evolution of inverse strain energy $1/\Psi$ over density $\rho$ given by Equation \ref{eq:integration_interp} for several starting volume fractions $v_0$.}
    \label{fig:f1-curves_log}
\end{figure}

A 20-point dataset is fed to the D$^2$NN model, half of them generated with $r=1$ (Figure \ref{fig:f2-r1_dataset_PINN} and the other half with $r=1$ (Figure \ref{fig:f3-r4_dataset_PINN}). Each set of 10 points is denoted with letters $a$-$j$ in Figure \ref{fig:f1-curves_log} showcasing curves for different initial volume fractions along which data points P1, P2, P3 and P4 are located. Whereas P1 and P2 are within the training dataset, P3 and P4 are external to it. The points $d$ and $e$ are also used for an example of interpolation within a curve.

\begin{figure}[H]
    \centering
    \includegraphics[width=0.95\linewidth]{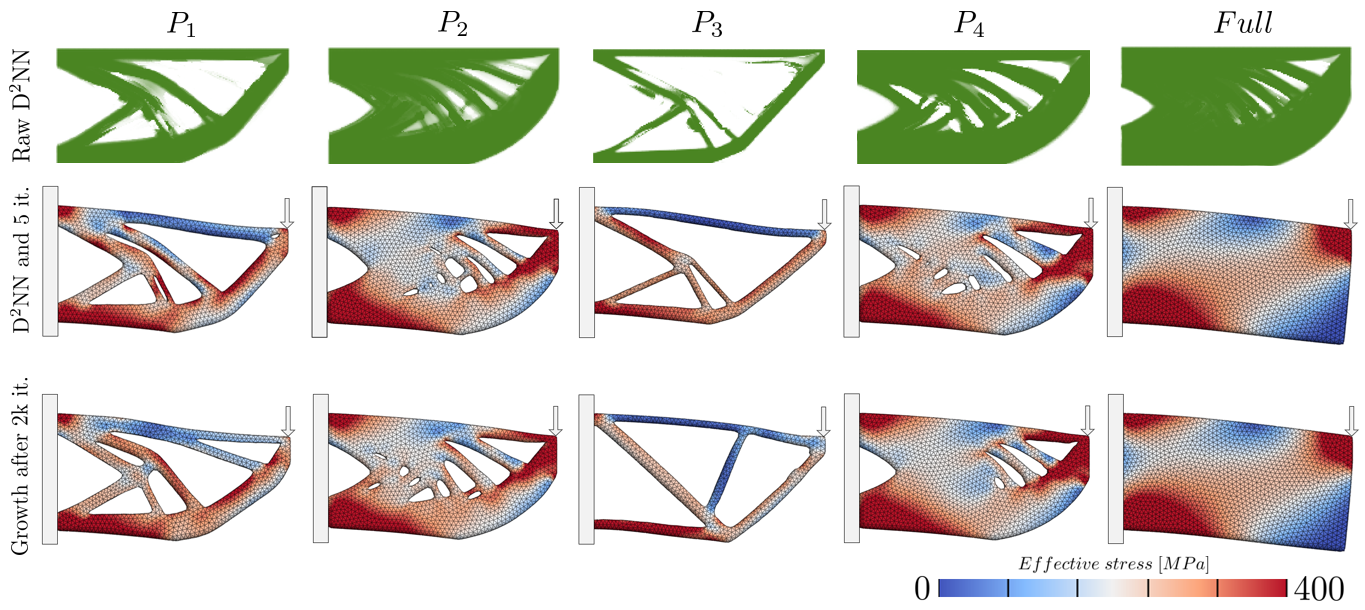}
    \caption{\centering Topology optimization of a cantilever beam for points $P_1$-$P_4$ in Figure \ref{fig:f1-curves_log} and the full material block (columns), via D$^2$NN exclusively (top row), D$^2$NN after 5 initial iterations of Equation \ref{eq:iterative_scheme} (middle row) and 2,000 iterations (bottom row). Blue-red scale indicates stress values (MPa). Erosion filtering radius $r = 4$.}
    \label{fig:f5-PINN_interpolation_fig}
\end{figure}

These images depict interpolation points A-J in Figure \ref{fig:f1-curves_log}, serving as database to obtain intermediate topologies $P_1$-$P_4$ via Equation \ref{eq:iterative_scheme}, D$^2$NN or both - see Figure \ref{fig:f5-PINN_interpolation_fig}. The physical loss in the D$^2$NN scheme is given by the distance (in compliance and volume) between the interpolation point's and the target's respective curves.

The results of training the D$^2$NN model with Figures \ref{fig:f2-r1_dataset_PINN} and \ref{fig:f3-r4_dataset_PINN} and testing it on unseen data points $P_1$-$P_4$ in Figure \ref{fig:f1-curves_log} is displayed in Figure \ref{fig:f5-PINN_interpolation_fig}, using FEBIO\textcopyright\ as a visualization tool and representing von Mises stress on stainless steel prototypes (AISI 316L, $E$ = 193 GPa, $\nu$ = 0.33).

Figure \ref{fig:f5-PINN_interpolation_fig} exhibits great column-wise similarities, meaning the D$^2$NN-based surrogate (500 training epochs, top row) reaches very close results to the vanilla scheme (2,000 iterations, bottom row), saving a lot of computation time since, once trained, the D$^2$NN-driven option is about 25 times faster. 

If the D$^2$NN scheme is only run for its 5 first iterations (middle row) rather than doing so from the initial curve's volume $v_0$, the results when applying Equation \ref{eq:iterative_scheme} on the predicted $\rho$ are even closer to the ground truth, with almost exactly matching stress distributions - especially in the full configuration.

Of course, this resemblance is never exact, since the dataset (A-J) is made up from points in different curves than the targets $P_1$-$P_4$ - see Figures \ref{fig:f2-r1_dataset_PINN} and \ref{fig:f3-r4_dataset_PINN} for examples with filtering radii $r = 1$ and $r = 4$, respectively, excluding the full volume block as the convergence point of all curves. Thus, this methodology could be interpreted as a quick way to obtain mechanically equivalent designs to the target - in terms of compliance, volume and minimum thickness - without actually having to compute all the needed iterations along their unexplored curves.

\begin{figure}[H]
    \centering
    \includegraphics[width=1\linewidth]{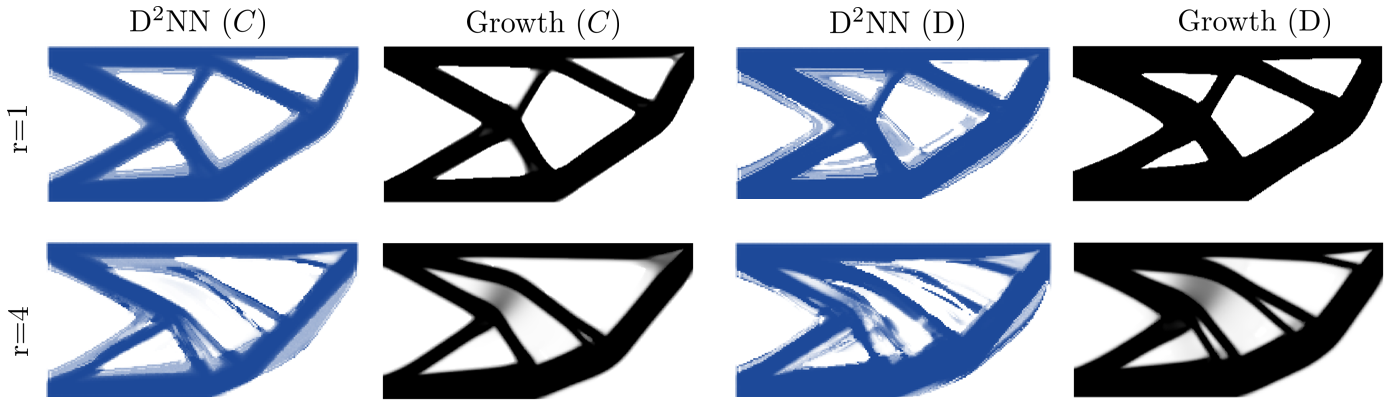}
    \caption{\centering Data point ground truth (logarithmic growth with Equation \ref{eq:iterative_scheme}, in black) and D$^2$NN prediction (in blue) along the training dataset curve A-E in Figure \ref{fig:f1-curves_log} (points C and D) for filtering radii $r = 1$ (top row) and $r = 4$ (bottom row).}
    \label{fig:f4-PINN_vs_growth}
\end{figure}

This tool can also be employed to predict points along the same A-E and F-J curves, quickening dataset generation for training and so further accelerating the whole TO approximation process. See Figure \ref{fig:f4-PINN_vs_growth} for some such examples: predictions are now much more accurate, since interpolation is done along the same curve for datasets and targets. This disrupts SIMP's initial homogeneous assumption, accelerating convergence.

\begin{figure}[H]
    \centering
    \includegraphics[width=1\linewidth]{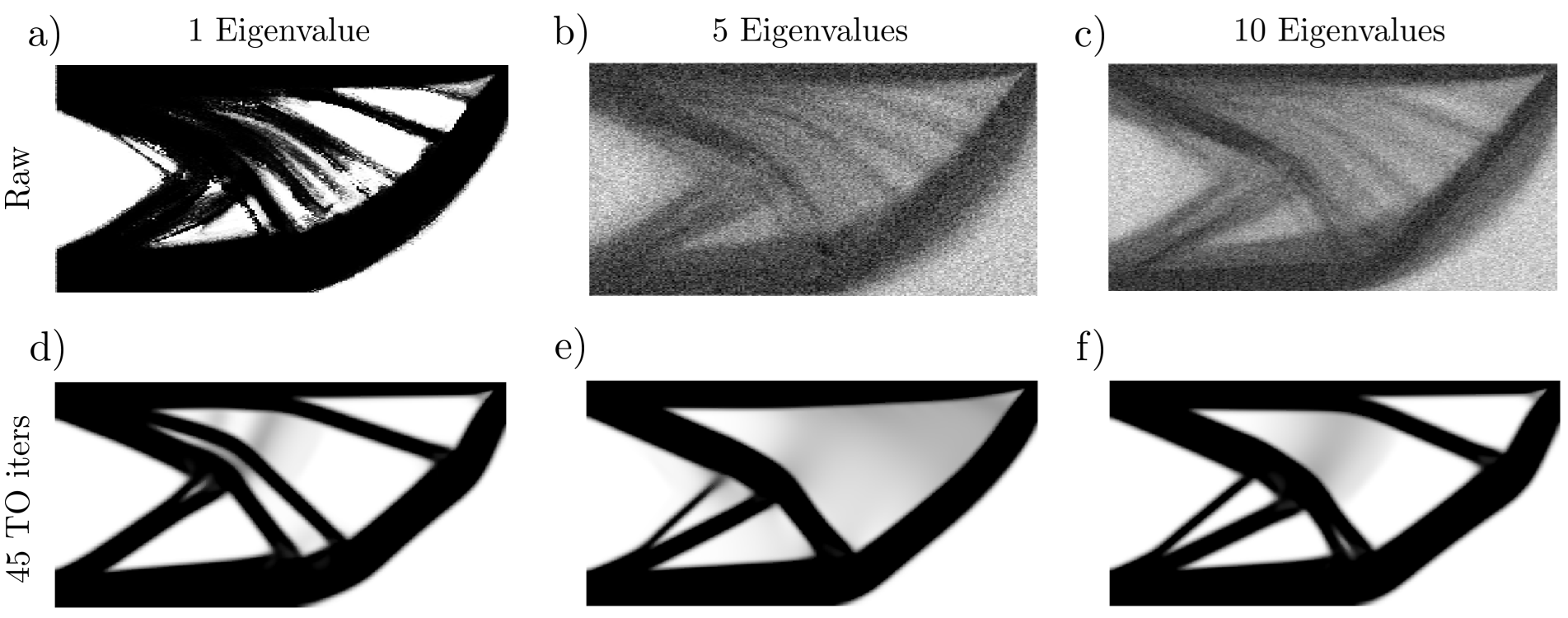}
    \caption{\centering D$^2$NN forecasts (top row) and 45 logarithmic iterations (Equation \ref{eq:iterative_scheme}, bottom row) with the first 1 (left column), 5 (middle column) and 10 eigenvalues (right column).}
    \label{fig:f6-PINN_P1_r4_with_eigen}
\end{figure}

If D$^2$NN training is performed upon a reduced order dataset (e.g. via Singular Value Decomposition), the topology's resolution can be tuned to produce slightly different results with virtually equivalent mechanical properties at a lower computational cost. See Figure \ref{fig:f6-PINN_P1_r4_with_eigen} for some examples of D$^2$NN-predicted (raw) and 45 iterations of Equation \ref{eq:iterative_scheme} for reduced eigenvalue sets, suitable for Principal Component Analysis.\\

Columns in Figures \ref{fig:f5-PINN_interpolation_fig} and \ref{fig:f6-PINN_P1_r4_with_eigen} are instances of diverse topologies for the same volume fraction $\rho$ and strain energy $\Psi$ in the $\Lambda$ space representing a sole $\Phi(\eta,\rho)$ point. As explained in Figure \ref{fig:Mickey}, this is but a consequence of the non-convex, highly non-linear solution space in TO, in which many initial points can lead to virtually equivalent designs in the $\Phi$ and/or $\Lambda$ spaces, highlighting the need for a richer $\Gamma(\eta,\rho,\alpha)$ space able to tell topologies $\alpha$ apart.






\section{Conclusions and future research lines}

This article has introduced some shape and mechanical issues affecting most common topological optimization schemes, namely capillary struts which, although mathematically correct (following the minimization process), do not efficiently distribute neither material nor mechanical loading (strain/stress), making them prone to stress concentration and thus vulnerable to fracture, among other problems; apart from difficult to manufacture depending on the chosen technique.

Two novel approaches have been offered as ways to alleviate such drawbacks. The first one consists of a density filter which effectively eliminates those useless thin ribs while redistributing an equal or even smaller volume more wisely, which in turn improves mechanical performance. The second proposal consists of an analytical expression linking the two main design variables in the SIMP method (compliance and volume fraction) in an experimentally-consistent manner which in turn provides different “families” evolving in shape and topology. Importantly, these two ingredients can be combined and allow for interpolation, as seen in various examples, and introduced in a Machine Learning surrogate scheme. 

These tools have been put to the test under various boundary conditions and loadings, proving they can create a versatile design space $\Lambda$ with adjustable compliance $c$, volume fraction $v$ and minimum thickness (via $r$); all of that within reasonable computational cost. This is deemed of great interest for practical reasons, since this article explains the full methodology for simultaneous 2D shape and topology inverse design with tuneable twig width (minimum thickness). Future improvements include a 3D filter generalization, an extension to probabilistic loads (varying in position, module and frequency) and further enrichment with embedded damage criteria (compliance element-wise penalization \cite{IrastorzaValera2025prob}).

%
%
%
%


\section*{Funding}

\includegraphics[width=0.15\linewidth]{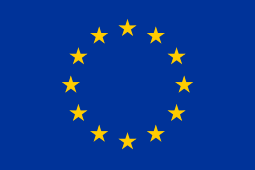}\includegraphics[width=0.33\linewidth]{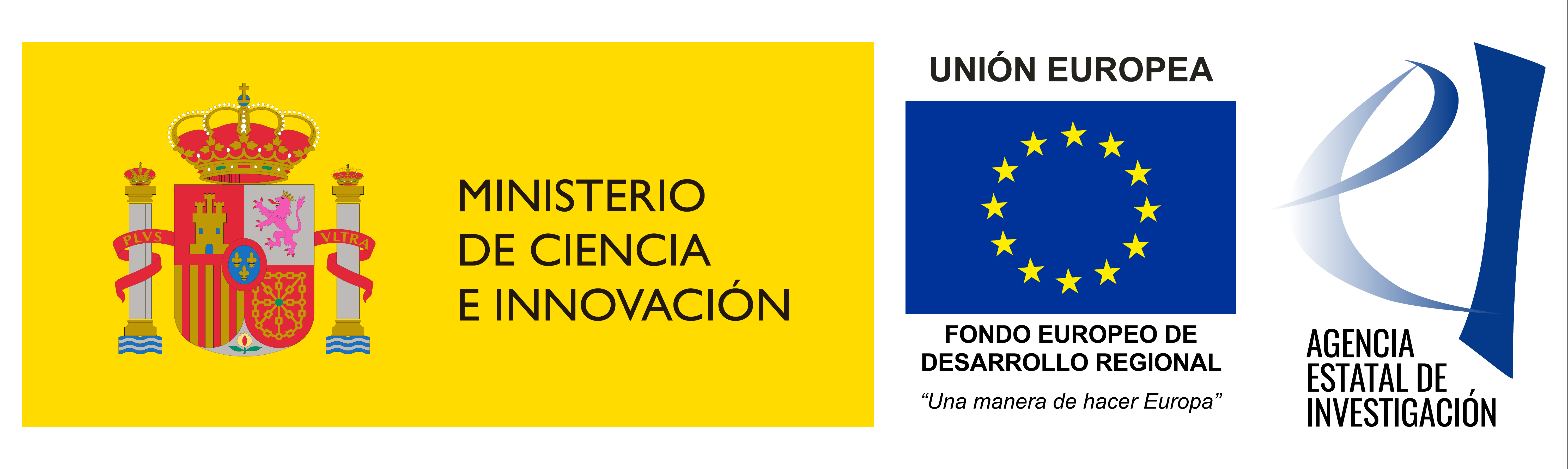}

This project has received funding from the European Union’s Horizon 2020 Marie Sk\l{}odowska-Curie Actions - Innovative European Training Networks under grant agreement No 956401, as well as from the Spanish Ministry of Science and Innovation and the State Research Agency (AEI) PID2021-126051OB-C43.\\



\section*{Conflict of interest}

The authors declare that the research was conducted in the absence of any commercial or financial relationships that could be construed as a potential conflict of interest.







\bibliographystyle{plain}

\bibliography{article}







\end{document}